\documentclass[11pt,leqno]{amsart}%
\usepackage{amsmath}
\usepackage{mathptmx}
\usepackage{amsfonts}
\usepackage{amssymb}
\usepackage{graphicx}%
\usepackage{hyperref}
\usepackage{color}
\usepackage{cancel}

\providecommand{\U}[1]{\protect\rule{.1in}{.1in}}
\newtheorem*{acknowledgement}{Acknowledgement}

\newtheorem{conjecture}{Conjecture}
\newtheorem{definition}{Definition}
\newtheorem{lemma}{Lemma}
\newtheorem{proposition}{Proposition}

\newtheorem{theorem}{Theorem}

\numberwithin{equation}{section}

\title[Critical metrics on 4-manifolds]{Critical metrics of the total scalar curvature functional on 4-manifolds}
\author{A. Barros, B. Leandro \& E. Ribeiro Jr.}
\address[A. Barros] {Universidade Federal do Cear\'a - UFC, Departamento  de Matem\'atica, Campus do Pici, Av. Humberto Monte, Bloco 914,
60455-760, Fortaleza / CE, Brazil.}\email{abbarros@mat.ufc.br}
\address[B. Leandro]{Universidade de Bras\'ilia - UNB, Departamento de Matem\'atica, 70910-900, Bras\'ilia / DF, Brazil.}\email{B.L.Neto@mat.unb.br}
\address[E. Ribeiro Jr]{Universidade Federal do Cear\'a - UFC, Departamento  de Matem\'atica, Campus do Pici, Av. Humberto Monte, Bloco 914,
60455-760, Fortaleza / CE, Brazil.}\email{ernani@mat.ufc.br}

\allowdisplaybreaks
\numberwithin{equation}{section}
\numberwithin{theorem}{section}
\thanks{A. Barros was partially supported by grant from CNPq/Brazil}
\thanks{B. Leandro was partially supported by grant from  CNPq/Brazil Proc. 149896/2012-3}
\thanks{E. Ribeiro Jr acknowledges partial support by CNPq/Brazil and Funcap/Brazil}
\keywords{Total scalar curvature functional, critical point equation,
Einstein metric.}\subjclass[2000]{Primary 53C25, 53C20, 53C21; Secondary 53C65}
\date{May 5, 2015}

\begin{document}
\newcommand{\spacing}[1]{\renewcommand{\baselinestretch}{#1}\large\normalsize}
\spacing{1.}

\begin{abstract}
The purpose of this paper is to investigate the critical points of the total scalar curvature functional restricted to space of metrics with constant scalar curvature of unitary volume, for simplicity CPE metrics.  It was conjectured in $1980$'s that every CPE metric must be Einstein.  We prove that a $4$-dimensional CPE metric with harmonic tensor $W^+$ must be isometric to a round sphere $\Bbb{S}^4.$
\end{abstract}

\maketitle

\section{Introduction}\label{introduction}

A fundamental problem in differential geometry is to find Riemannian metrics on a given manifold that provides constant curvature. An useful tool in this direction is to analyze the critical points of the total scalar curvature functional. More precisely, let $M^n$ be a compact oriented smooth manifold and $\mathcal{M}$ the
set of smooth Riemannian structures on $M^n$ of unitary volume. Given a
metric $g\in \mathcal{M}$ we define the total scalar curvature, or Einstein-Hilbert
functional $S:\mathcal{M}\to\Bbb{R}$ by
\begin{equation}
\label{functional} \displaystyle{S(g)=\int_{M}R_{g}dM_{g},}
\end{equation}
where $R_g$ and $dM_{g}$ stand, respectively, for the scalar
curvature of $M^n$ and the volume form. Einstein and Hilbert have proved that
the critical points of the functional $S$ are Einstein, for more details see Theorem 4.21 in \cite{besse}. The Einstein-Hilbert functional  restricted to a given conformal class is just the Yamabe functional, whose critical points are constant scalar curvature metrics in that class.

A result obtained by combining results due to  Aubin \cite{aubin}, Schoen \cite{schoen}, Trudinger \cite{trudinger} and Yamabe \cite{yamabe} gives the existence of a constant scalar curvature metric in every conformal class of Riemannian metrics on a compact manifold $M^n.$ Therefore, it is interesting to consider the set $$\mathcal{C}=\{g\in \mathcal{M};\,R_g\,\, \hbox{is \,constant}\}.$$ In \cite{koiso}, Koiso showed that, under  generic condition, $\mathcal{C}$ is an infinite dimensional manifold (cf. Theorem 4.44 in \cite{besse} p. 127).

It is well-known that the formal $L^{2}$-adjoint $\mathfrak{L}_{g}^{*}$ of the linearization $\mathfrak{L}_{g}$ of the scalar curvature
operator at $g$ is given by
\begin{equation}
\label{eqlinear}
\mathfrak{L}_{g}^{*}(f)=-(\Delta_{g}f)g+Hess_{g} f-f Ric_{g},
\end{equation} where $f$ is a smooth function on $M^n.$ Moreover, at any given metric $g$ in the space of the metrics with constant scalar curvature the map $\mathfrak{L}_{g}^{*}$  defined on $C^{\infty}$ to $\mathcal{M}$ is an over determined elliptic operator (cf. \cite{anderson}). Formally the Euler-Lagrangian equation of Hilbert-Einstein action restricted to $\mathcal{C}$ may be written as the following critical point equation

\begin{equation}
\label{eqf}
\mathfrak{L}_{g}^{*}(f)=Ric-\frac{R}{n}g,
\end{equation} where $f$ is a smooth function on $M^n.$ From this, it is easy to check that (\ref{eqf}) becomes

\begin{equation}\label{eqfund}
Ric-\frac{R}{n}g=Hess f -\big(Ric-\frac{R}{n-1} g\big)f,
\end{equation}where $Ric,$ $R$ and $Hess $ stand, respectively, for the Ricci tensor, the scalar curvature and the Hessian form on $M^n.$ Obviously Einstein metrics are recovered when $f=0.$ Moreover, the existence of a non constant solution is only known in the round sphere for some height function. It has been conjectured that the critical points of the total scalar curvature functional $S$ restricted to $\mathcal{C}$ are Einstein. For more details, we refer the reader to \cite{besse} p. 128.

\begin{definition}
A CPE metric is a 3-tuple $(M^n,\,g,\,f),$ where $(M^n,g)$ is a compact oriented Riemannian manifold of dimension $n\geq3$
with constant scalar curvature while $f$ is a smooth potential satisfying equation (\ref{eqfund}).
\end{definition}

Note that, computing the trace in (\ref{eqfund}), we obtain
\begin{equation}
\label{tracefund} \Delta f+\frac{R}{(n-1)}f=0.
\end{equation}Whence $f$ is an eigenfunction of the Laplacian and then $R$ is positive.

The conjecture proposed in  \cite{besse} in the middle of $1980$'s may be restated in terms of CPE (see e.g. \cite{br}, \cite{hwang2010}, \cite{Hwang1}, \cite{hwang}, \cite{benedito} and \cite{jiewei}). More exactly, the quoted authors proposed the following conjecture.

\begin{conjecture}[1980's, \cite{besse}]
\label{conj}
 A CPE metric is always Einstein.
\end{conjecture}

In the last years many authors have been tried to settle up this conjecture, but only partial results were achieved.  Among them, we detach the next ones: Lafontaine \cite{lafontaine} in 1983 proved it under locally conformally flat assumption and $Ker\,\mathfrak{L}_{g}^{*}(f)\neq 0$ and in 2011  Chang, Hwang and Yun avoided the condition on $Ker\,\mathfrak{L}_{g}^{*}(f),$  (see  Corollary 1.3  in \cite{CHY}); moreover, they also settled it  for metrics with parallel Ricci tensor, see \cite{CHY1}.  In 2000, Hwang \cite{Hwang1} was able to obtain the conjecture provided $f\ge-1$ and in 2003 \cite{hwang} he gave a topological answer showing that a three-dimensional compact manifold with null second homology group such that $Ker\, \mathfrak{L}_{g}^{*}(f) \neq 0$ must be diffeomorphic to a round sphere $\Bbb{S}^{3}$ and in 2010  Chang, Hwang and Yun \cite{hwang2010} showed that under these last conditions the conjecture is true. In 2014, Chang, Hwang and Yun \cite {CHY} proved the conjecture under harmonic curvature assumption.

It is well-known that 4-dimensional compact Riemannian manifolds have special behavior. In large part, this is because the bundle of $2$-forms on a 4-dimensional compact oriented Riemannian manifold can be invariantly decomposed as a direct sum; many relevant facts may be found in \cite{besse} and \cite{scorpan}.  For instance, on an oriented Riemannian manifold $(M^4,\,g),$ the Weyl curvature tensor $W$ is an endomorphism of the bundle of 2-forms $\Lambda^2=\Lambda^{2}_{+}\oplus\Lambda^{2}_{-}$ such that  $$W = W^+\oplus W^-,$$ where $W^\pm : \Lambda^{2}_\pm \longrightarrow \Lambda^{2}_\pm$ are called of the self-dual and anti-self-dual parts of $W.$ Half conformally flat metrics are also known as self-dual or anti-self-dual if $W^{-}$ or $W^{+}=0,$ respectively. Barros and Ribeiro Jr \cite{br} showed that Conjecture \ref{conj} is also true for $4$-dimensional half conformally flat manifolds. We highlight that $\Bbb{CP}^2$ endowed with Fubini-Study metric shows that the half-conformally flat condition is weaker than locally conformally flat condition in dimension 4. While Qing and Yuan \cite{jiewei} obtained a positive answer for Bach-flat manifolds in any dimension.

Proceeding, viewing $W^+$ as a tensor of type $(0,4),$ we say that the tensor $W^+$ is harmonic if $\delta W^+ = 0,$ where $\delta$ is the formal divergence defined for any
$(0,4)$-tensor $T$ by $$\delta T(X_1,X_2,X_3) =
trace_{g}\{(Y,Z)\mapsto\nabla_{Y}T(Z,X_1,X_2,X_3)\},$$ where $g$ is
the metric of $M^4.$ It is important to highlight that in dimension 4 we have $$|\delta W|^2=|\delta W^{+}|^2+|\delta W^{-}|^2.$$ Moreover, it should be emphasized that every 4-dimensional Einstein manifold has harmonic tensor $W^+$ (cf. 16.65 in \cite{besse}, see also Lemma 6.14 in \cite{handbook}). We also recall that every oriented 4-dimensional manifolds with $W^+$ harmonic satisfies the following relation
\begin{equation}
\label{weitzenbock}
\Delta |W^{+}|^2=2|\nabla W^{+}|^2+R|W^{+}|^2- 36 \det W^{+}.
\end{equation} For details see, for instance, Proposition 16.73 in \cite{besse}. Therefore, it is natural to ask which geometric implications has the assumption of the harmonicity of the tensor $W^+$ on a CPE metric.

Inspired by the historical development on the study of CPE conjecture, in this paper, we prove that, in dimension 4, the assumption that $(M^4 ,\,g)$ is locally conformally flat, considered in \cite{lafontaine}  and \cite{CHY}, as well as half conformally flat considered in \cite{br}, can be replaced by the weaker condition that $M^4$ has harmonic tensor $W^{+}.$ More precisely, we have the following result.

\begin{theorem}
\label{thmA} Conjecture \ref{conj} is true for $4$-dimensional manifolds with harmonic tensor $W^+.$
\end{theorem}

The main tools in the proof of the above result is to analyze the behavior of the level sets of the potential function $f$ which defines a CPE metric combined with some pointwise arguments. Obviously if we change the condition $\delta W^{+} = 0$ by the conditon $\delta W^{-} = 0$ the conclusion of Theorem \ref{thmA} is exactly the same. Furthermore, one should be emphasized that there is no relationship between Bach-flat condition, considered in  \cite{jiewei}, and the condition that $M^4$ has harmonic tensor $W^{+}.$

This article is organized as follows. In Section 2, we review some classical tensors which we shall use here. Moreover, we outline some useful informations about 4-dimensional manifolds. In Section 3, we prove the main result.

\section{Background}

Throughout this section we recall some basic tensors and informations that will be useful in the proof of our main result.  First of all, for operators $S,T:\mathcal{H} \to \mathcal{H}$ defined over an $n$-dimensional Hilbert space $\mathcal{H}$ the Hilbert-Schmidt inner product is defined according to
\begin{equation}
\langle S,T \rangle =\rm tr\big(ST^{\star}\big),
\label{inner}
\end{equation}where $\rm tr$ and $\star$ denote, respectively, the trace and the adjoint operation. Moreover, if $I$ denotes the identity operator on $\mathcal{H}$ the traceless operator of $T$ is given by
 \begin{equation}
 \label{eqtr1}
\mathring{T}=T - \frac{\rm tr T}{n}I.
 \end{equation} In particular the norm of $\mathring{T}$ satisfies
 \begin{equation}
 \label{eqnorm1}
 \mid \mathring{T} \mid^2 = \mid T \mid^2 -\frac{(\rm tr T)^2}{n}.
 \end{equation}

Now, we recall that for a Riemannian manifold $(M^{n},\,g),$ $n\geq 3,$ the Weyl tensor $W$ is defined by the following decomposition formula
\begin{eqnarray}
\label{weyl}
R_{ijkl}&=&W_{ijkl}+\frac{1}{n-2}\big(R_{ik}g_{jl}+R_{jl}g_{ik}-R_{il}g_{jk}-R_{jk}g_{il}\big)\nonumber\\&&-\frac{R}{(n-1)(n-2)}\big(g_{jl}g_{ik}-g_{il}g_{jk}\big),
\end{eqnarray} where $R_{ijkl}$ stands for the Riemannian curvature operator. Moreover, the Cotton tensor $C$ is given according  to
\begin{equation}
\label{cotton}
\displaystyle{C_{ijk}=\nabla_{i}R_{jk}-\nabla_{j}R_{ik}-\frac{1}{2(n-1)}\big(\nabla_{i}R g_{jk}-\nabla_{j}R g_{ik}).}
\end{equation}
These two tensors are related as follows
\begin{equation}
\label{cottonwyel} \displaystyle{C_{ijk}=-\frac{(n-2)}{(n-3)}\nabla_{l}W_{ijkl},}
\end{equation}provided $n\ge 4.$
Finally, the Schouten tensor  $A$ is given by
\begin{equation}
\label{schouten}
A_{ij}=\frac{1}{n-2}\big(R_{ij}-\frac{R}{2(n-1)}g_{ij}\big).
\end{equation} For more details about these tensors we address  to \cite{besse}.

In what follows $M^4$ will denote an oriented 4-dimensional manifold and $g$ is a Riemannian metric on $M^4.$ As it was previously mentioned 4-manifolds are fairly special.  For instance, following the notations used in \cite{handbook} (see also \cite{VIA} p. 46), given any local orthogonal frame $\{e_{1}, e_{2}, e_{3}, e_{4}\}$ on open set of $M^4$ with dual basis $\{e^{1}, e^{2}, e^{3}, e^{4}\},$ there exists a unique bundle morphism $\ast$ called {\it Hodge star} (acting on bivectors), such that $$\ast(e^{1}\wedge e^{2})=e^{3}\wedge e^{4}.$$ This implies that $\ast$ is an involution, i.e. $\ast^{2}=Id.$ In particular, this implies that the bundle of $2$-forms on a 4-dimensional oriented Riemannian manifold can be invariantly decomposed as a direct sum $\Lambda^2=\Lambda^{2}_{+}\oplus\Lambda^{2}_{-}.$ This allows us to conclude that the Weyl tensor $W$ is an endomorphism of $\Lambda^2=\Lambda^{+} \oplus \Lambda^{-} $ such that
\begin{equation}
\label{decW}
W = W^+\oplus W^-.
\end{equation} We recall that $dim_{\Bbb{R}}(\Lambda^2)=6$ and $dim_{\Bbb{R}}(\Lambda^{\pm})=3.$ Also, it is well-known that
\begin{equation}
\label{6A}
\Lambda^{+}=span\Big\{\frac{e^{1}\wedge e^{2}+e^{3}\wedge e^{4}}{\sqrt{2}},\,\frac{e^{1}\wedge e^{3}+e^{4}\wedge e^{2}}{\sqrt{2}},\,\frac{e^{3}\wedge e^{2}+e^{4}\wedge e^{1}}{\sqrt{2}}\Big\}
\end{equation}
 and
\begin{equation}
\label{6B}
 \Lambda^{-}=span\Big\{\frac{e^{1}\wedge e^{2}-e^{3}\wedge e^{4}}{\sqrt{2}},\,\frac{e^{1}\wedge e^{3}-e^{4}\wedge e^{2}}{\sqrt{2}},\,\frac{e^{3}\wedge e^{2}-e^{4}\wedge e^{1}}{\sqrt{2}}\Big\}.
\end{equation} From this, the bundles $\Lambda^{+}$ and $\Lambda^{-}$ carry natural orientations such that the bases (\ref{6A}) and (\ref{6B}) are positive-oriented. Furthermore, if $\mathcal{R}$ denotes the curvature of $M^4$ we get a matrix

\begin{equation}
\mathcal{R}=
\left(
  \begin{array}{c|c}
    \\
W^{+} +\frac{R}{12}Id & \mathring{Ric} \\ [0.4cm]\hline\\

    \mathring{Ric}^{\star} & W^{-}+\frac{R}{12}Id  \\[0.4cm]
  \end{array}
\right),
\end{equation}
where $\mathring{Ric}:\Lambda^{-}\to \Lambda^{+}$ stands for the Ricci traceless operator of $M^4.$

Recalling that the Weyl tensor is trace-free on any pair of indices we have
\begin{equation}
\label{W+}
W_{\,p\,q\,r\,s}^{+}=\frac{1}{2}\big(W_{p\,q\,r\,s}+W_{p\,q\,\overline{r}\,\overline{s}}\big),
\end{equation} where $(\overline{r}\,\overline{s}),$ for instance, stands for the dual of $(r\,s),$ that is, $(r\,s\,\overline{r}\,\overline{s})=\sigma(1234)$ for some even permutation $\sigma$ in the set $\{1,2,3,4\}$ (cf. Equation 6.17, p. 466 in \cite{handbook}). In particular, we have $$W_{1234}^{+}=\frac{1}{2}\big(W_{1234}+W_{1212}\big).$$ For more details we refer to \cite{handbook} and \cite{VIA}.

\section{Proof of the main result}

In order to set the stage for the proof to follow let us recall an useful result obtained in \cite{br} for any dimension.

\begin{lemma}[\cite{br}]
\label{lem21}
Let $\big(M^n,\,g,\,f)$ be a CPE metric. Then:
\begin{eqnarray*}
(f+1)C_{ijk}&=&W_{ijks}\nabla^{s}f-\frac{R}{(n-2)}\big(\nabla_{j}f g_{ik} - \nabla_{i}f g_{jk}\big)+\frac{(n-1)}{(n-2)}\big(R_{ik}\nabla_{j}f-R_{jk}\nabla_{i}f\big)\\&&-\frac{1}{(n-2)}\big(R_{is}\nabla^{s}f g_{jk}-R_{js}\nabla^{s}fg_{ik}\big).
\end{eqnarray*}
\end{lemma} Next, following the notation used in \cite{br} we define the tensor $T$ as follows

\begin{eqnarray}
\label{tensorT}
T_{ijk}&=&\frac{(n-1)}{(n-2)}\big(R_{ik}\nabla_{j}f-R_{jk}\nabla_{i}f\big)-\frac{1}{(n-2)}\big(R_{is}\nabla^{s}f g_{jk}-R_{js}\nabla^{s}fg_{ik}\big)\nonumber\\&&-\frac{R}{(n-2)}\big(\nabla_{j}f g_{ik} - \nabla_{i}f g_{jk}\big).
\end{eqnarray}

Taking into account this definition we deduce from Lemma \ref{lem21} that
\begin{equation}
\label{eq59}
(f+1)C_{ijk}=W_{ijks}\nabla^{s}f+T_{ijk}.
\end{equation}

Next, we need the following results by Hwang et al. \cite{Hwang1,hwang}.

 \begin{proposition}[\cite{Hwang1,hwang}]
\label{prop1a}
 Let $(M^n,\,g,\,f)$ be a CPE metric with $f$ non-constant. Then, the set $\{x\in M^{n}: f(x) = -1\}$ has measure zero.
\end{proposition}

We now denote by $Crit(f)$ the set $\{x\in M^{n} ; \nabla f(x) = 0\}.$ With this setting, we have the following proposition.

\begin{proposition}[\cite{CHY}]\label{prop45}
Let $(M^n,\,g,\,f),$ be a CPE metric. Then $Crit(f)$ has zero $n$-dimensional measure.
\end{proposition}

Now we are ready to prove Theorem \ref{thmA}.

\subsection{Proof of Theorem \ref{thmA}}

\begin{proof}

To begin with, we shall compute the value of $\delta W^{+}$ on a 4-dimensional CPE metric. In fact, since a CPE metric has constant scalar curvature, Equation (\ref{cotton}) becomes
 $$C_{klj}=\nabla_{k}R_{lj} - \nabla_{l}R_{kj}.$$ So, as an immediate consequence of (\ref{cottonwyel}) we have

\begin{equation}
\label{lemk1}
2\delta W_{jkl}^{+}=\frac{1}{2}\big(\nabla_{k}R_{jl}-\nabla_{l}R_{jk}\big)+\frac{1}{2}\big(\nabla_{\overline{k}}R_{j\overline{l}}-\nabla_{\overline{l}}R_{j\overline{k}}\big).
\end{equation} From Lemma \ref{lem21} and (\ref{tensorT}) we already know that
\begin{equation}
\label{CWT}
(f+1)C_{ijk}=W_{ijks}\nabla^{s}f+T_{ijk}.
\end{equation} This combined with (\ref{lemk1}) gives
\begin{eqnarray}
4(1+f)\delta W^{+}_{jkl}&=& (1+f)\big(\nabla_{k}R_{jl}-\nabla_{l}R_{jk}\big)+(1+f)\big(\nabla_{\overline{k}}R_{j\overline{l}}
-\nabla_{\overline{l}}R_{j\overline{k}}\big)\nonumber\\
&=&(1+f)(C_{klj}+C_{\overline{k}\,\overline{l}\,j})\nonumber\\
&=&\left[W_{kljs}\nabla^{s}f + W_{\overline{k}\,\overline{l}\,j\,s}\nabla^{s}f + T_{lkj}+T_{\overline{l}\,\overline{k}\,j}\right]
\end{eqnarray}
Now, using that $\delta W^{+}=0$ we infer
\begin{eqnarray}
-W_{kljs}\nabla^{s}f - W_{\overline{k}\,\overline{l}\,j\,s}\nabla^{s}f= T_{klj}+T_{\overline{k}\,\overline{l}\,j}.
\end{eqnarray}
Hence we have
\begin{eqnarray}
0=-\left(W_{kljs}\nabla^{s}f + W_{\overline{k}\,\overline{l}\,j\,s}\nabla^{s}f\right)\nabla^{j}f= \left(T_{klj}+T_{\overline{k}\,\overline{l}\,j}\right)\nabla^{j}f,
\end{eqnarray} so that
\begin{eqnarray}
\label{T1}
\left(T_{ijk}+T_{\overline{i}\,\overline{j}\,k}\right)\nabla^{k}f=0.
\end{eqnarray}

On the other hand, Equation (\ref{tensorT}) allows us to deduce
\begin{eqnarray}
T_{ijk}\nabla^{k}f&=&\frac{3}{2}(R_{ik}\nabla_{j}f -R_{jk}\nabla_{i}f)\nabla^{k}f - \frac{1}{2}(R_{is}\nabla^{s}fg_{jk}-R_{js}\nabla^{s}fg_{ik})\nabla^{k}f\nonumber\\
&&-\frac{R}{2}(\nabla_{j}fg_{ik}-\nabla_{i}fg_{jk})\nabla^{k}f\nonumber\\&=&(R_{ik}\nabla_{j}f -R_{jk}\nabla_{i}f)\nabla^{k}f.
\end{eqnarray}Then, from (\ref{T1}) we arrive at
\begin{eqnarray}
\label{kl}
(R_{ik}\nabla_{j}f -R_{jk}\nabla_{i}f)\nabla^{k}f + (R_{\overline{i}\,k}\nabla_{\overline{j}}f -R_{\overline{j}\,k}\nabla_{\overline{i}}f)\nabla^{k}f=0.
\end{eqnarray}

In the sequel, we consider an orthonormal frame $\{e_{1},e_{2},e_{3},e_{4}\}$ diagonalizing $Ric$ at a point $q,$ such that $\nabla f (q)\neq 0,$ with associated eigenvalues $\lambda_{k},\,(k=1,\ldots,4),$ respectively. It is important to highlight that the regular points of $M^4,$ denoted by $\{p\in M^n:\nabla f(p)\ne 0\},$ is dense in $M^4.$  Otherwise, $f$ must be constant in an open set of $M^4;$ in fact, from  (\ref{tracefund}) this constant is zero, but from standard theory of nodal sets $f$ can not vanish in an open set, see for instance \cite{cheng}. From now on, up to explicit mention, we restrict our attention to regular points. So, a straightforward computation using (\ref{kl}) gives the following useful system
\begin{eqnarray}\label{sytem}
\left\{
\begin{array}{lcc}
(\lambda_{1}-\lambda_{2})\nabla_{1}f\nabla_{2}f+(\lambda_{3}-\lambda_{4})\nabla_{3}f\nabla_{4}f=0,\\
(\lambda_{1}-\lambda_{3})\nabla_{1}f\nabla_{3}f+(\lambda_{4}-\lambda_{2})\nabla_{4}f\nabla_{2}f=0,\\
(\lambda_{1}-\lambda_{4})\nabla_{1}f\nabla_{4}f+(\lambda_{2}-\lambda_{3})\nabla_{2}f\nabla_{3}f=0.
\end{array}
\right.
\end{eqnarray}

We now claim that $\nabla f,$ whenever nonzero, is an eigenvector for $Ric.$ In fact, taking into account that $\nabla f(p) \neq0$ we have that at least one of the $(\nabla_{j}f)\neq0$, $1\leq j\leq4$. If this occurs for exactly one of them, then $\nabla f=(\nabla_{j}f)e_j$ for some $j$, which gives that $Ric(\nabla f)=\lambda_j \nabla f.$ On the other hand, if we have $(\nabla_{j}f)\neq0$  for two directions, without loss of generality we can suppose that $\nabla_{1}f\neq0$, $\nabla_{2}f\neq0$, $\nabla_{3}f=0$ and $\nabla_{4}f=0.$ Then, from (\ref{sytem}) we have $\lambda_{1}=\lambda_{2}=\lambda$. In a such case we have $\nabla f=(\nabla_{1}f)e_{1} +(\nabla_{2}f)e_{2}.$ From this, we obtain
\begin{eqnarray}
Ric(\nabla f) &=&Ric((\nabla_{1}f)e_{1} +(\nabla_{2}f)e_{2})=(\nabla_{1}f)Ric(e_{1}) +(\nabla_{2}f)Ric(e_{2})\nonumber\\
&=&(\nabla_{1}f)\lambda_{1}e_{1} +(\nabla_{2}f)\lambda_{2}e_{2}=\lambda\nabla f.
\end{eqnarray} Next, the case $(\nabla_{j}f)\neq0$ for three directions is analogous. Now, it remains to analyze the case $(\nabla_{j}f)\neq0$ for $j=1,2,3,4.$ In this case we use once more (\ref{sytem}) to obtain
\begin{eqnarray}
(\lambda_{1}-\lambda_{2})^{2}(\nabla_{1}f\nabla_{2}f)^{2}&+&(\lambda_{3}-\lambda_{4})^{2}(\nabla_{3}f\nabla_{4}f)^{2}\nonumber\\
+(\lambda_{1}-\lambda_{3})^{2}(\nabla_{1}f\nabla_{3}f)^{2}&+&(\lambda_{4}-\lambda_{2})^{2}(\nabla_{4}f\nabla_{2}f)^{2}\nonumber\\
+(\lambda_{1}-\lambda_{4})^{2}(\nabla_{1}f\nabla_{4}f)^{2}&+&(\lambda_{2}-\lambda_{3})^{2}(\nabla_{2}f\nabla_{3}f)^{2}=0.
\end{eqnarray} Therefore, $\lambda_{1}=\lambda_{2}=\lambda_{3}=\lambda_{4}.$ From here it follows that $\nabla f$ is an eigenvector for $Ric,$ which proves our claim.

Proceeding, we study the level sets of the potential function $f$ which defines a CPE metric. To that end we denote $\Sigma_{c}=\{p\in M^n : f(p)=c\}.$ At regular points the vector field $e_{1}=\frac{\nabla f}{\mid \nabla f \mid}$ is normal to $\Sigma_{c}$ and $\{e_{2},e_{3},e_{4}\}$ is  an orthonormal frame on $\Sigma_{c}.$ With this notation in mind, since $Ric(\nabla f)=\lambda \nabla f$ and $\nabla_{e_{a}}f=\langle \nabla f, e_{a}\rangle=0$ for $a=\{2,3,4\},$ we immediately deduce from (\ref{tensorT}) that  $T_{abc}=0$ for $\{a,b,c\}=\{2,3,4\}.$ Whence, in $\Sigma_{c},$ (\ref{eq59}) becomes

\begin{equation}
\label{89}
(f+1)C_{abc}=W_{abcs}\nabla^{s}f.
\end{equation}

We notice that, for an arbitrary $Y\in T_{p}\Sigma_{c},$ we have $\mathring{Ric}(\nabla f, Y)=0.$ So, we can use the fundamental equation to infer

$$\nabla_{e_{a}}|\nabla f|^{2}=(f+1)\mathring{Ric}(\nabla f,e_{a})-\frac{Rf}{n(n-1)}g(\nabla f,e_{a})=0.$$ This implies that $|\nabla f|^{2}$ is constant in $\Sigma_{c}.$

On the other hand, it is well-known that the second fundamental form of $\Sigma$ is given by 

\begin{equation}
\label{eq123456}
h_{ab}=\langle \nabla_{e_{a}}e_{b}, e_{1}\rangle=-\frac{f_{ab}}{|\nabla f|},
\end{equation} where $\{a,b\}=\{2,3,4\}.$ Moreover, its mean curvature is $$H_{\Sigma}=\frac{1}{|\nabla f|^{2}}\big(f_{11}-\Delta f\big).$$ We then combine the fundamental equation with (\ref{eq123456}) to conclude that near the point $p$ the second fundamental form is given by $$h_{ab}=\mu g_{ab},$$ where $\mu$ is equal to the mean curvature. But, since $|\nabla f|^{2}$ is constant in $\Sigma_{c}$ follows that $H$ is constant on $\Sigma,$ so is $\mu$ (see also \cite{br}). Then, we may use Codazzi's equation to infer
\begin{eqnarray}
R_{abc1}=\nabla^{\Sigma}_{b}h_{ac}-\nabla^{\Sigma}_{a}h_{bc}=0
\end{eqnarray} at $p$ when $\{a,b,c\}=\{2,3,4\}.$ In particular, from (\ref{weyl}) we infer $W_{abc1}=0$ and this jointly with (\ref{89}) gives $$(f+1)C_{abc}=0.$$ We then use Proposition \ref{prop1a} to conclude that $C_{abc}=0$ when $\{a,b,c\}=\{2,3,4\}.$

Next, we already know from (\ref{lemk1}) that
\begin{eqnarray}
4\delta W^{+}_{jkl}=C_{klj}+C_{\bar{k}\bar{l}j}.
\end{eqnarray} Therefore, since $\delta W^{+}=0$ and $C_{abc}=0$ we get $C_{\bar{a}\bar{b}c}=0.$ More precisely, we have
\begin{eqnarray}
0=C_{\bar{1}\bar{2}c}=C_{34c},\quad\quad
0=C_{\bar{1}\bar{3}c}=-C_{24c}\quad\mbox{and}\quad
0=C_{\bar{2}\bar{3}c}=C_{14c}.
\end{eqnarray} Similarly, we can use (\ref{T1}) jointly with (\ref{eq59}) to conclude that $C_{ij1}=0$ when $\{i,j\}=\{1,2,3,4\}.$ This allows us to conclude that $C_{ijk}=0$ when $\{i,j,k\}=\{1,2,3,4\}.$ Whence, we have from  (\ref{cottonwyel}) that $\delta W=0$ in $M^4.$ Since $M^4$ has constant scalar curvature, then $M^4$ has harmonic curvature. Finally, it suffices to invoke Theorem 1.2 in \cite{CHY} to conclude that $M^4$ is isometric to a round sphere $\Bbb{S}^{4}.$ In addition, we also conclude that $f$ is a height function on $\Bbb{S}^{4}.$ 

This finishes the proof of Theorem \ref{thmA}.
\end{proof}

\begin{acknowledgement}
The authors want to thank the referees for their careful reading and helpful suggestions. Moreover, the third author wish to express his gratitude for the excellent support during his stay at Department of Mathematics - Lehigh University, where part of this paper was carried out.
\end{acknowledgement}

\end{document}